\documentclass[12pt]{article}

\usepackage{amsfonts}
\usepackage[centertags]{amsmath}
\usepackage{tikz-cd}
\usepackage{enumerate}
\usepackage[pdftex]{pict2e}
\usepackage{xcolor}

\usepackage{graphicx}
\usepackage{import}
\usepackage{epstopdf}

\def\squarebox#1{\hbox to #1{\hfill\vbox to #1{\vfill}}}
\newcommand{\qed}{\hspace*{\fill}
\vbox{\hrule\hbox{\vrule\squarebox{.667em}\vrule}\hrule}\smallskip}
\newtheorem{teorema}{Theorem}[section]
\newtheorem{lema}[teorema]{Lemma}

\newtheorem{proposicao}[teorema]{Proposition}

\newenvironment{prova}{\noindent {\bf Proof:}}{\hfill $\qed $ \newline}

\newcommand{\vol}{\mbox{{\rm vol}}}

\begin{document}

\title{Fully geometric definitions of \\
orientation and determinants}

\author{
Mauro Patr\~{a}o\footnote{Departament of Mathematics,
University of Brasília, Brazil. \textit{mpatrao@mat.unb.br }.
}
}

\maketitle

\begin{abstract}
We present fully geometric definitions of orientation and determinants and show they coincide with the algebraic definitions. This allows us to provide an approach to determinants in the spirit of what is presented in the article A Geometric Approach to Determinants, previously published by the Monthly, but making it rigorous and complete.
\end{abstract}
  


The relations between volumes and determinants are well known, but the attempts to introduce determinants via volumes are not very successful. In this article, we provide an approach to determinants in the spirit of what is presented in \cite{hannah}, but making it rigorous and complete. We present a fully geometric definition of determinants and show that it coincides with the algebraic definition using two simple geometric properties satisfied by the signed volumes of parallelepipeds. One of the major challenges is to provide a fully geometric definition of orientation which enables us to derive rigorously its main properties.

\section{Preliminaries}

A fully geometric definition means a rigorous definition based only in geometric concepts such as scalar product, euclidean distance, translations, rotations, trigonometric functions, and continuity.
Given noncollinear unitary vectors $u$ and $w$, we can define
\[
 \widehat{u} = \frac{u - (w \cdot u)w}{|u - (w \cdot u)w|}
\]
and also
\begin{eqnarray*}
 R_{u,w}^tv & = & v - (\widehat{u} \cdot v)\widehat{u} - (w \cdot v)w \\
 & & +(\cos(t)(\widehat{u} \cdot v)-\sin(t)(w \cdot v))\widehat{u} \\
 & & +(\sin(t)(\widehat{u} \cdot v)+\cos(t)(w \cdot v))w
\end{eqnarray*}
We have that $(t,u,v,w) \mapsto R_{u,w}^tv$ is a continuous map and that $v \mapsto R_{u,w}^tv$ is a rotation in the plane generated by $u$ and $w$ such that
\[
 R_{u,w}^0v = v
 \qquad \mbox{e} \qquad
 R_{u,w}^{T_{u,w}}u = w
\]
where
\[
 T_{u,w} = \frac{\widehat{u} \cdot u}{|\widehat{u} \cdot u|}\mbox{acos}(w \cdot u)
\]
and $\mbox{acos}$ denotes the arccosine.

\begin{center}
\begin{picture}(200,180)
\put(0,0){\includegraphics[scale=.3]{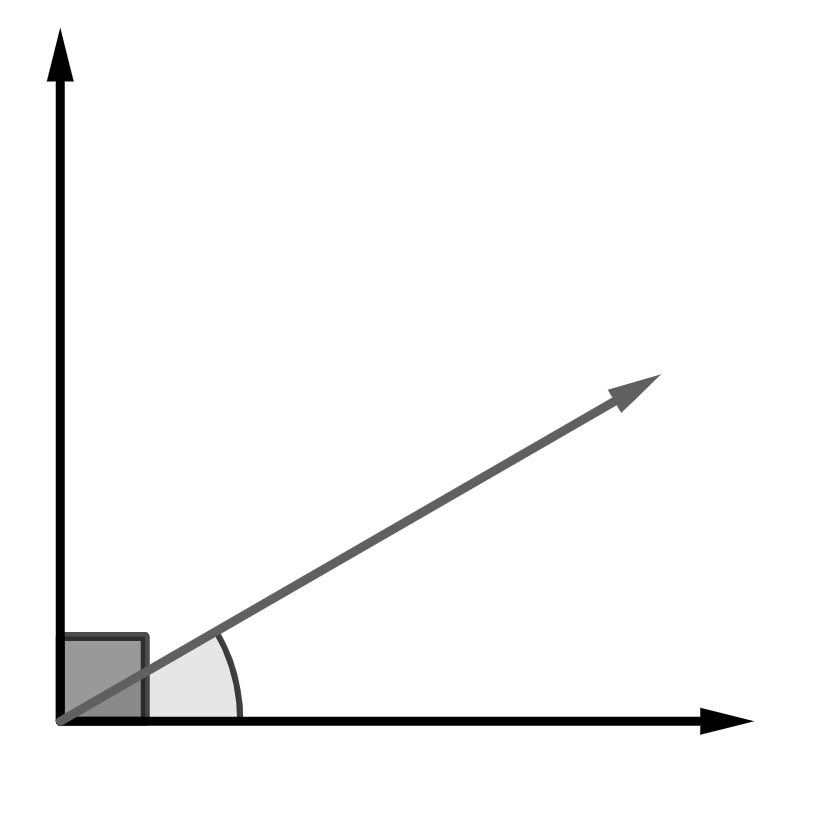}}
\put(170,8){$ w = R_{u,w}^{T_{u,w}}u$}
\put(0,170){$\widehat{u}$}
\put(100,80){$u$}
\put(60,30){$\mathrm{acos}{({u}\cdot{w})}$}
\end{picture}

\end{center}

We need the following results in order to provide a fully geometric definition of orientation.

\begin{lema}\label{lemma-bases-path}
Given $(u_1,\ldots,u_n)$ and $(w_1,\ldots,w_n)$, with $n \geq 2$, two sufficiently closed ordered orthonormal bases, there exists a continuous path $t \mapsto (b_1(t),\ldots,b_n(t))$ of orthonormal ordered bases such that
\begin{equation}
 b_i(0) = u_i, \qquad b_i(1) = w_i
\end{equation}
for every $i \in \{1,\ldots,n\}$.
\end{lema}

\begin{prova}
We proceed by induction on $n$, the dimension of the vector space.

When $n=2$, we have two possibilities. First, if $u_2 = w_2$, then $u_1 = w_1$, since the two bases are sufficiently closed each other. In this case, we can define $b_i(t) = u_i = w_i $ for every $i \in \{1,2\}$ and every $t \in [0,1]$. Second, if $u_2 \neq w_2$, then $u_2$ and $w_2$ are noncollinear, since $u_2 \neq -w_2$, because the two bases are sufficiently closed each other. In this case, we can define $b_i(t) = R_{u_2,w_2}^{t}u_i$ for every $i \in \{1,2\}$ and every $t \in [0,T_2]$, where $T_2 = T_{u_2,w_2}$. It follows that $b_i(0) = u_i$ for every $i \in \{1,2\}$ and that $b_2(T_2) = R_{u_2,w_2}^{T_2}u_2 = w_2$. Since the two bases are sufficiently closed each other, it follows that $b_1(T_2) = R_{u_2,w_2}^{T_2}u_1 = w_1$. The claim follows after normalizing the parameter $t$.

When $n>2$, we also have two possibilities. First, if $u_n = w_n$, we have that $(u_1,\ldots,u_{n-1})$ and $(w_1,\ldots,w_{n-1})$ are two sufficiently closed ordered orthonormal bases of the subspace perpendicular to $u_n = w_n$, which is a vector space of dimension $n-1$. By the induction hypothesis, there exists a continuous path $t \mapsto (b_1(t),\ldots,b_{n-1}(t))$ of orthonormal ordered bases of the subspace perpendicular to $u_n = w_n$ such that $b_i(0) = u_i$ and $b_i(1) = w_i$ for every $i \in \{1,\ldots,n-1\}$. In this case, we can define $b_n(t) = u_n = w_n$ for every $t \in [0,1]$. Second, if $u_n \neq w_n$, then $u_n$ and $w_n$ are noncollinear, since $u_n \neq -w_n$, because the two bases are sufficiently closed each other. In this case, we can define $\widehat{b}_i(t) = R_{u_n,w_n}^{t}u_i$ for every $i \in \{1,\ldots,n\}$ and every $t \in [0,T_n]$, where $T_n = T_{u_n,w_n}$. It follows that $\widehat{b}_i(0) = u_i$ for every $i \in \{1,\ldots,n\}$ and that $\widehat{b}_n(T_n) = R_{u_n,w_n}^{T_n}u_n = w_n$. We have that $(R_{u_n,w_n}^{T_n}u_1,\ldots,R_{u_n,w_n}^{T_n}u_{n-1})$ and $(w_1,\ldots,w_{n-1})$ are two sufficiently closed ordered orthonormal bases of the subspace perpendicular to $R_{u_n,w_n}^{T_n}u_n = w_n$, which is a vector space of dimension $n-1$. By the induction hypothesis, there exists a continuous path $t \mapsto (b_1(t),\ldots,b_{n-1}(t))$ of orthonormal ordered bases of the subspace perpendicular to $R_{u_n,w_n}^{T_n}u_n = w_n$ such that $b_i(0) = R_{u_n,w_n}^{T_n}u_i$ and $b_i(1) = w_i$ for every $i \in \{1,\ldots,n-1\}$. The claim follows defining $b_n(t) = R_{u_n,w_n}^{T_n}u_n = w_n$, connecting continuously $\widehat{b}_i(t)$ and $b_i(t)$ for every $i \in \{1,\ldots,n\}$ and normalizing the parameter $t$.
\end{prova}

\begin{proposicao}
Given $(u_1,\ldots,u_n)$ and $(w_1,\ldots,w_n)$, with $n \geq 3$, two sufficiently closed ordered orthonormal bases such that $u_n, w_n \neq \pm e_n$, there exists a continuous path $t \mapsto (c_1(t),\ldots,c_n(t))$ of orthonormal ordered bases such that
\begin{equation}
 c_i(0) = u_i, \qquad c_i(1) = w_i, \qquad c_n(t) \neq \pm e_n
\end{equation}
for every $i \in \{1,\ldots,n\}$ and every $t \in [0,1]$.
\end{proposicao}

\begin{prova}
We have two possibilities. First, if $u_n = w_n$, we have that $(u_1,\ldots,u_{n-1})$ and $(w_1,\ldots,w_{n-1})$ are two sufficiently closed ordered orthonormal bases of the subspace perpendicular to $u_n = w_n$, which is a vector space of dimension $n-1 \geq 2$. By Lemma \ref{lemma-bases-path}, there exists a continuous path $t \mapsto (b_1(t),\ldots,b_{n-1}(t))$ of orthonormal ordered bases of the subspace perpendicular to $u_n = w_n$ such that $b_i(0) = u_i$ and $b_i(1) = w_i$ for every $i \in \{1,\ldots,n-1\}$. In this case, we can define $c_i(t) = b_i(t)$ for every $i \in \{1,\ldots,n-1\}$ and $c_n(t) = u_n = w_n \neq \pm e_n$ for every $t \in [0,1]$.

Second, if $u_n \neq w_n$, there exists $v$ distinct from but closed to $u_n$ and $w_n$ such that $\pm e_n$ is not in the plane generated by $u_n$ and $v$, nor in the plane generated by $w_n$ and $v$. In this case, we can define $\widehat{b}_i(t) = R_{u_n,v}^{t}u_i$ for every $i \in \{1,\ldots,n\}$ and every $t \in [0,\widehat{T}]$, where $\widehat{T} = T_{u_n,v}$. It follows that $\widehat{b}_i(0) = u_i$ for every $i \in \{1,\ldots,n\}$, that $\widehat{b}_n(\widehat{T}) = R_{u_n,v}^{\widehat{T}}u_n = v$, and that $\widehat{b}_n(t) \neq \pm e_n$ for every $t \in [0,\widehat{T}]$. Now, we can define $\widetilde{b}_i(t) = R_{v,w_n}^{t}\widehat{b}_i(\widehat{T})$ for every $i \in \{1,\ldots,n\}$ and every $t \in [0,\widetilde{T}]$, where $\widetilde{T} = T_{v,w_n}$. It follows that $\widetilde{b}_i(0) = \widehat{b}_i(\widehat{T})$ for every $i \in \{1,\ldots,n\}$, that $\widetilde{b}_n(\widetilde{T}) = R_{v,w_n}^{\widetilde{T}}\widehat{b}_n(\widehat{T}) = R_{v,w_n}^{\widetilde{T}}v = w_n$, and that $\widetilde{b}_n(t) \neq \pm e_n$ for every $t \in [0,\widetilde{T}]$. We have that $(\widetilde{b}_1(\widetilde{T}),\ldots,\widetilde{b}_{n-1}(\widetilde{T}))$ and $(w_1,\ldots,w_{n-1})$ are two sufficiently closed ordered orthonormal bases of the subspace perpendicular to $\widetilde{b}_n(\widetilde{T}) = w_n$, which is a vector space of dimension $n-1 \geq 2$. By Lemma \ref{lemma-bases-path}, there exists a continuous path $t \mapsto (b_1(t),\ldots,b_{n-1}(t))$ of orthonormal ordered bases of the subspace perpendicular to $\widetilde{b}_n(\widetilde{T}) = w_n$ such that $b_i(0) = \widetilde{b}_i(\widetilde{T})$ and $b_i(1) = w_i$ for every $i \in \{1,\ldots,n-1\}$. The claim follows defining $b_n(t) = \widetilde{b}_n(\widetilde{T}) = w_n$, and defining $c_i(t)$ by connecting continuously $\widehat{b}_i(t)$, $\widetilde{b}_i(t)$, and $b_i(t)$ for every $i \in \{1,\ldots,n\}$ and normalizing the parameter $t$.
\end{prova}

\section{Orientation}

Given an ordered orthonormal basis $(e_1,\ldots,e_n)$, we define inductively the orientation of a given ordered orthonormal basis $(u_1,\ldots,u_n)$ with respect to $(e_1,\ldots,e_n)$. In the unidimensional case, we have that $u_1=\pm e_1$, and hence we can define the orientation of $u_1$ with respect to $e_1$ as follows
\[
S_1(u_1) =
 \left\{
 \begin{array}{rl}
  1, & \quad u_1 = e_1 \\
  -1, & \quad u_1 = -e_1
 \end{array}
 \right.
\]
In the bidimensional case, we have that
\begin{eqnarray*}
 c_1(t)
 & = & +(\cos(t)(e_1 \cdot u_1)-\sin(t)(e_2 \cdot u_1))e_1 \\
 & & +(\sin(t)(e_1 \cdot u_1)+\cos(t)(e_2 \cdot u_1))e_2
 \\
 c_2(t)
 & = & +(\cos(t)(e_1 \cdot u_2)-\sin(t)(e_2 \cdot u_2))e_1 \\
 & & +(\sin(t)(e_1 \cdot u_2)+\cos(t)(e_2 \cdot u_2))e_2
\end{eqnarray*}
are such that $(c_1(t),c_2(t))$ is an ordered orthonormal basis for every $t$ and such that
\[
 c_1(0) = u_1, \quad c_2(0) = u_2, \quad c_1(T_2) = \pm e_1, \quad c_2(T_2) = e_2
\]
where $T_2 = ((e_1 \cdot u_2)/|e_1 \cdot u_2|)\mbox{acos}(e_2 \cdot u_2)$. Hence we can define the orientation of $(u_1,u_2)$ with respect to $(e_1,e_2)$ as follows
\[
S_2(u_1,u_2) =
 \left\{
 \begin{array}{rl}
  1, & \quad c_1(T_2) = e_1 \\
  -1, & \quad c_1(T_2) = -e_1
 \end{array}
 \right.
\]
which is clearly a continuous map.

\begin{center}

\begin{picture}(320,180)
\put(0,0){\includegraphics[scale=.3]{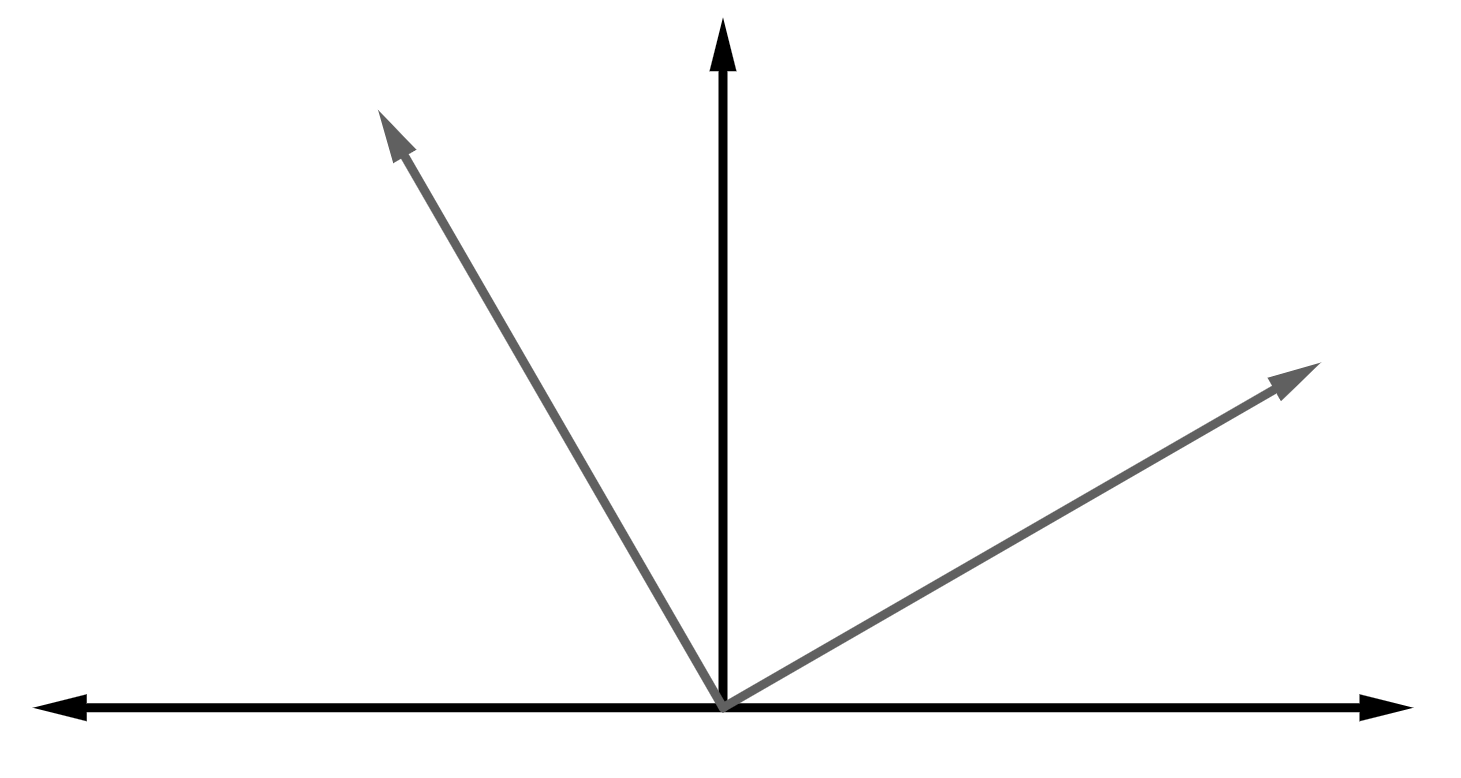}}
\put(320,0){$e_1$}
\put(170,170){$e_2 = c_2(T_2)$}
\put(30,150){$c_1(0)=u_1$}
\put(-30,0){$c_1(T_2)=-e_1$}
\put(300,100){$u_2 = c_2(0)$}

\end{picture}

\end{center}

For $n \geq 3$ and when $u_n \neq \pm e_n$, we can define inductively the orientation of $(u_1,\ldots,u_n)$ with respect to $(e_1,\ldots,e_n)$ as follows
\[
 S_n(u_1,\ldots,u_n) = S_{n-1}\left(R_{u_n,e_n}^{T_n}u_1,\ldots,R_{u_n,e_n}^{T_n}u_{n-1}\right)
\]
where $T_n = T_{u_n,e_n}$ and $\left(R_{u_n,e_n}^{T_n}u_1,\ldots,R_{u_n,e_n}^{T_n}u_{n-1}\right)$ is an ordered orthonormal basis of the space generated by the basis $(e_1,\ldots,e_{n-1})$, since $R_{u_n,e_n}^{T_n}u_n = e_n$.

\begin{center}

\begin{picture}(320,250)
\put(0,0){\includegraphics[scale=.8]{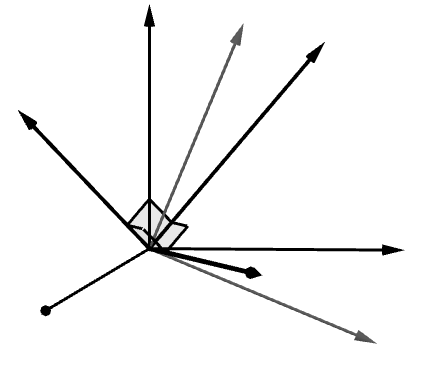}}
\put(200,190){$u_1$}
\put(100,220){$e_2$}
\put(0,160){$u_3$}
\put(-20,20){$R_{u_3,e_3}^{T_3}u_3=e_3$}
\put(150,210){$R_{u_3,e_3}^{T_3}u_1$}
\put(230,10){$R_{u_3,e_3}^{T_3}u_2$}
\put(160,50){$u_2$}
\put(250,60){$e_1$}

\end{picture}

\end{center}

Assuming inductively that $S_{n-1}$ is a continuous map, we have that the map $S_n(u_1,\ldots,u_n)$ is continuous at every ordered orthonormal basis $(u_1,\ldots,u_n)$, where $u_n \neq \pm e_n$. Besides this, given an ordered orthonormal basis $(v_1,\ldots,v_n)$, where $v_n = \pm e_n$, and given two other ordered orthonormal bases $(u_1,\ldots,u_n)$ and $(w_1,\ldots,w_n)$ sufficiently close to the first one, with $u_n, w_n \neq \pm e_n$, there exists a continuous path $t \mapsto (c_1(t),\ldots,c_n(t))$ of orthonormal ordered bases such that
\[
 c_i(0) = u_i, \qquad c_i(1) = w_i, \qquad c_n(t) \neq \pm e_n
\]
for every $i \in \{1,\ldots,n\}$ and every $t \in [0,1]$. Thus $t \mapsto S_n(c_1(t),\ldots,c_n(t))$ is a continuous map from $[0,1]$ to $\{\pm 1\}$, and hence it is constant. Therefore $S_n(u_1,\ldots,u_n)=S_n(w_1,\ldots,w_n)$, showing that $S_n$ is constant in a neighborhood of $(v_1,\ldots,v_n)$, where $v_n = \pm e_n$, and hence we can define
\[
 S_n(v_1,\ldots,v_n) = S_n(u_1,\ldots,u_n)
\]
for any $(u_1,\ldots,u_n)$ sufficiently close to $(v_1,\ldots,v_n)$ with $u_n \neq \pm e_n$. Thus we have that $S_n(u_1,\ldots,u_n)$ is locally constant and hence continuous at every  ordered orthonormal basis $(u_1,\ldots,u_n)$.

\begin{proposicao}\label{lemmasignal}
\begin{equation}
 S_n(-u_1,u_2,\ldots,u_n) = -S_n(u_1,u_2,\ldots,u_n)
\end{equation}
\end{proposicao}

\begin{prova}
 We proceed by induction in $n$. For $n=2$, the claim is immediate. For $n\geq3$, assuming that the claim is valid for $n-1$, we have that
 \begin{eqnarray*}
  S_n(-u_1,u_2,\ldots,u_n)
  &=& S_{n-1}\left(R_{u_n,e_n}^{T_n}-u_1,R_{u_n,e_n}^{T_n}u_2,\ldots,R_{u_n,e_n}^{T_n}u_{n-1}\right) \\
  &=& S_{n-1}\left(-R_{u_n,e_n}^{T_n}u_1,R_{u_n,e_n}^{T_n}u_2,\ldots,R_{u_n,e_n}^{T_n}u_{n-1}\right) \\
  &=& -S_{n-1}\left(R_{u_n,e_n}^{T_n}u_1,R_{u_n,e_n}^{T_n}u_2,\ldots,R_{u_n,e_n}^{T_n}u_{n-1}\right) \\
  &=& -S_n(u_1,u_2,\ldots,u_n)
 \end{eqnarray*}
when $u_n \neq \pm e_n$, and the result follows in the general case, since $S_n$ is locally constant.
\end{prova}

For an ordered basis $(v_1,\ldots,v_n)$, not necessarily orthonormal, we define
\[
 G(v_1,\ldots,v_n) = (u_1,\ldots,u_n)
\]
given by the Gram-Schmidt orthornomalization process, where $u_1 = v_1/|v_1|$ and inductively
\[
 u_i = \frac{v_i - \sum_{j=1}^{i-1}(v_i\cdot u_j)u_j}{\left|v_i - \sum_{j=1}^{i-1}(v_i\cdot u_j)u_j\right|}
\]
We define the orientation of $(v_1,\ldots,v_n)$ with respect to $(e_1,\ldots,e_n)$ as follows
\[
 S(v_1,\ldots,v_n) = S_n(G(v_1,\ldots,v_n))
\]
which is locally constant, since $S_n$ is locally constant and $G$ is continuous. We need the following bidimensional lemma, whose proof is immediate.

\begin{lema}\label{lemma2d}
Given $v_i$ and $v_j$ noncollinear, the exists a continuous path $t \mapsto \left(c_i(t),c_j(t)\right)$ of ordered bases of the space generated by $\{v_i, v_j\}$ such that
\begin{equation}
 \left(c_i(0),c_j(0)\right) = \left(v_i, v_j\right), \qquad \left(c_i(1),c_j(1)\right) = \left(-v_j, v_i\right)
\end{equation}
\end{lema}

\vspace{0.0cm}

\begin{center}

\begin{picture}(250,210)
\put(0,0){\includegraphics[scale=.2,angle=45]{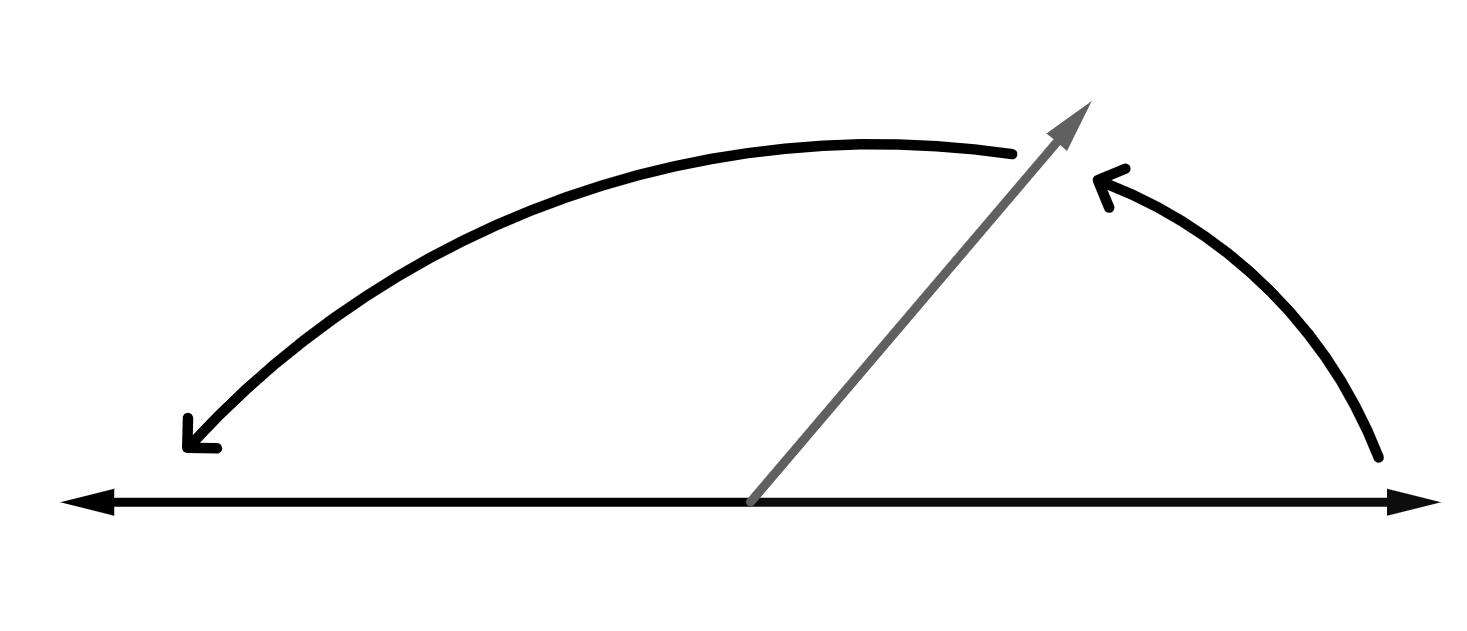}}
\put(150,180){$c_j(t)$}
\put(60,170){$c_i(0) = c_j(1)$}
\put(50,100){$c_i(t)$}
\put(-20,20){$c_i(1)=-c_j(0)$}
\put(210,170){$c_j(0)$}

\end{picture}

\end{center}

Next result is the main property of the orientation of ordered bases necessary in order to obtain the fundamental property of determinants.

\begin{proposicao}
\begin{equation}
 S(v_1,\ldots,v_{i-1},tv_i,v_{i+1},\ldots,v_n) = \frac{t}{|t|}S(v_1,\ldots,v_n)
\end{equation}
\end{proposicao}

\begin{prova}
We have that
\[
 G(tv_1,v_2,\ldots,v_n) = \left(\frac{t}{|t|}u_1,u_2,\ldots,u_n\right)
\]
if
\[
 G(v_1,\ldots,v_n) = \left(u_1,\ldots,u_n\right)
\]
which implies that
\begin{eqnarray*}
  S(tv_1,v_2,\ldots,v_n)
  &=& S_n\left(\frac{t}{|t|}u_1,u_2,\ldots,u_n\right) \\
  &=& \frac{t}{|t|}S_n\left(u_1,\ldots,u_n\right) \\
  &=& \frac{t}{|t|}S(v_1,\ldots,v_n)
 \end{eqnarray*}
where we used Proposition \ref{lemmasignal} in the second equality.Therefore it follows that
\begin{eqnarray*}
  S(v_1,\ldots,v_{i-1},tv_i,v_{i+1},\ldots,v_n)
  &=& S(-tv_i,\ldots,v_{i-1},v_1,v_{i+1},\ldots,v_n) \\
  &=& \frac{t}{|t|}S(-v_i,\ldots,v_{i-1},v_1,v_{i+1},\ldots,v_n) \\
  &=& \frac{t}{|t|}S(v_1,\ldots,v_{i-1},v_i,v_{i+1},\ldots,v_n) \\
  &=& \frac{t}{|t|}S(v_1,\ldots,v_n)
 \end{eqnarray*}
where we used Lemma \ref{lemma2d} in the first and third equalities.
\end{prova}

\section{Determinants}

We now introduce a geometric definition of the determinant of an ordered set  $\{v_1,\ldots,v_n\}$ as a signed volume of the parallelepiped generated by it. The parallelepiped generated by $\{v_1,\ldots,v_n\}$ is given by the set
\[
P(v_1,\ldots,v_n) = \{t_1v_1 + \cdots + t_nv_n : t_1,\ldots,t_n \in [0,1)\}
\]
and its volume is given by
\[
\vol(v_1,\ldots,v_n) = \mu(P(v_1,\ldots,v_n))
\]
where $\mu$ is the unique measure on $V$ invariant by translations such that $\mu(P(e_1,\ldots,e_n)) = 1$. Note that, for convenience, some faces of the parallelepiped are not included, but this does not change its volume. The volume of the parallelepiped satisfies the following simple properties:

\begin{enumerate}
\item[\textbf{V1}] For all vectors $v_1,\ldots,v_n$ and all nonnegative $t$, we have that
\[
\vol(v_1,\ldots,tv_i,\ldots,v_n) = t\vol(v_1,\ldots,v_i,\ldots,v_n)
\]
\begin{center}
\begingroup%
  \makeatletter%
  \providecommand\color[2][]{%
    \errmessage{(Inkscape) Color is used for the text in Inkscape, but the package 'color.sty' is not loaded}%
    \renewcommand\color[2][]{}%
  }%
  \providecommand\transparent[1]{%
    \errmessage{(Inkscape) Transparency is used (non-zero) for the text in Inkscape, but the package 'transparent.sty' is not loaded}%
    \renewcommand\transparent[1]{}%
  }%
  \providecommand\rotatebox[2]{#2}%
  \ifx\svgwidth\undefined%
    \setlength{\unitlength}{179.34337812bp}%
    \ifx\svgscale\undefined%
      \relax%
    \else%
      \setlength{\unitlength}{\unitlength * \real{\svgscale}}%
    \fi%
  \else%
    \setlength{\unitlength}{\svgwidth}%
  \fi%
  \global\let\svgwidth\undefined%
  \global\let\svgscale\undefined%
  \makeatother%
  \begin{picture}(1,0.6858501)%
    \put(0,0){\includegraphics[width=\unitlength,page=1]{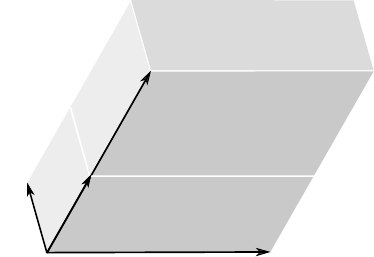}}%
    \put(0.74022622,0.0127674){\color[rgb]{0,0,0}\makebox(0,0)[lb]{\smash{$v_1$}}}%
    \put(0.25523891,0.15185833){\color[rgb]{0,0,0}\makebox(0,0)[lb]{\smash{$tv_2$}}}%
    \put(0.41134657,0.43752514){\color[rgb]{0,0,0}\makebox(0,0)[lb]{\smash{$v_2$}}}%
    \put(-0.00291946,0.22378153){\color[rgb]{0,0,0}\makebox(0,0)[lb]{\smash{$v_3$}}}%
  \end{picture}%
\endgroup%

\end{center}
\item[\textbf{V2}] For all vectors $v_1,\ldots,v_n$, we have that
\[
\vol(v_1,\ldots,v_i+v_j,\ldots,v_j,\ldots,v_n) = \vol(v_1,\ldots,v_i,\ldots,v_j,\ldots,v_n)
\]
\begin{center}
\begingroup%
  \makeatletter%
  \providecommand\color[2][]{%
    \errmessage{(Inkscape) Color is used for the text in Inkscape, but the package 'color.sty' is not loaded}%
    \renewcommand\color[2][]{}%
  }%
  \providecommand\transparent[1]{%
    \errmessage{(Inkscape) Transparency is used (non-zero) for the text in Inkscape, but the package 'transparent.sty' is not loaded}%
    \renewcommand\transparent[1]{}%
  }%
  \providecommand\rotatebox[2]{#2}%
  \ifx\svgwidth\undefined%
    \setlength{\unitlength}{290.16886351bp}%
    \ifx\svgscale\undefined%
      \relax%
    \else%
      \setlength{\unitlength}{\unitlength * \real{\svgscale}}%
    \fi%
  \else%
    \setlength{\unitlength}{\svgwidth}%
  \fi%
  \global\let\svgwidth\undefined%
  \global\let\svgscale\undefined%
  \makeatother%
  \begin{picture}(1,0.42127363)%
    \put(0,0){\includegraphics[width=\unitlength,page=1]{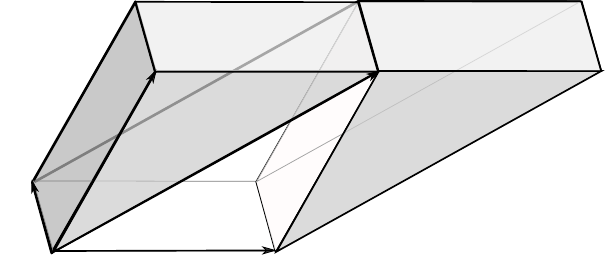}}%
    \put(0.38493019,0.02128986){\color[rgb]{0,0,0}\makebox(0,0)[lb]{\smash{$v_1$}}}%
    \put(0.25571711,0.26056061){\color[rgb]{0,0,0}\makebox(0,0)[lb]{\smash{$v_2$}}}%
    \put(-0.0020701,0.12784906){\color[rgb]{0,0,0}\makebox(0,0)[lb]{\smash{$v_3$}}}%
    \put(0.6265535,0.26531222){\color[rgb]{0,0,0}\makebox(0,0)[lb]{\smash{$v_1+v_2$}}}%
  \end{picture}%
\endgroup%

\end{center}
\end{enumerate}

These properties, which are very intuitive from the pictures, are proved rigorously in the appendix.
Now we can define the determinant of an ordered basis $(v_1,\ldots,v_n)$ as the product of its orientation by its volume
\[
 \det(v_1,\ldots,v_n) = S(v_1,\ldots,v_n)\mbox{vol}(v_1,\ldots,v_n)
\]
and, when $(v_1,\ldots,v_n)$ is not a base, its determinant is defined as being null. With this definition, we are able to obtain the fundamental properties of determinants.

\newpage

\begin{proposicao}
For every vectors $v_1,\ldots,v_n$ and every scalar $t$, we have that
\begin{enumerate}
\item[\textbf{D1})]
\begin{equation}
\det(v_1,\ldots,tv_i,\ldots,v_n) = t\det(v_1,\ldots,v_i,\ldots,v_n)
\end{equation}

\item[\textbf{D2})]
\begin{equation}
\det(v_1,\ldots,v_i+v_j,\ldots,v_j,\ldots,v_n) = \det(v_1,\ldots,v_i,\ldots,v_j,\ldots,v_n)
\end{equation}
\end{enumerate}
\end{proposicao}

\begin{prova}
When $(v_1,\ldots,v_n)$ is not a basis, we have that $(v_1,\ldots,tv_i,\ldots,v_n)$ as well $(v_1,\ldots,v_i+v_j,\ldots,v_j,\ldots,v_n)$ are not bases, and hence both side of the two above equations are null. When $(v_1,\ldots,v_n)$ is a basis, first we have that
\[
 \mbox{vol}(v_1,\ldots,-v_i,\ldots,v_n) = \mbox{vol}(v_1,\ldots,v_i,\ldots,v_n)
\]
since the parallelepiped generated by $(v_1,\ldots,v_i,\ldots,v_n)$ is the translation by $2v_i$ of the parallelepiped generated by $(v_1,\ldots,-v_i,\ldots,v_n)$. Hence it follows that
\begin{eqnarray*}
 \det(v_1,\ldots,tv_i,\ldots,v_n)
 &=& S(v_1,\ldots,tv_i,\ldots,v_n)\mbox{vol}(v_1,\ldots,tv_i,\ldots,v_n) \\
 &=& \frac{t}{|t|}S(v_1,\ldots,v_i,\ldots,v_n)|t|\mbox{vol}\left(v_1,\ldots,\frac{t}{|t|}v_i,\ldots,v_n\right) \\
 &=& tS(v_1,\ldots,v_i,\ldots,v_n)\mbox{vol}\left(v_1,\ldots,v_i,\ldots,v_n\right) \\
 &=& t\det(v_1,\ldots,v_i,\ldots,v_n)
\end{eqnarray*}
where we used property \textbf{V1} in the second equality and that $t/|t| = \pm 1$ in the the third equality. Besides this, we have that $t \mapsto c(t) = v_i+tv_j$ is a continuous path such that $c(0) = v_i$ and $c(1) = v_i+v_j$ which implies that
\[
 S(v_1,\ldots,v_i+v_j,\ldots,v_j,\ldots,v_n) = S(v_1,\ldots,v_i,\ldots,v_j,\ldots,v_n)
\]
and, by property \textbf{V2}, we have that
\[
 \mbox{vol}(v_1,\ldots,v_i+v_j,\ldots,v_j,\ldots,v_n) = \mbox{vol}(v_1,\ldots,v_i,\ldots,v_j,\ldots,v_n)
\]
Therefore, multiplying the two above equations, it follow that
\[
 \det(v_1,\ldots,v_i+v_j,\ldots,v_j,\ldots,v_n) = \det(v_1,\ldots,v_i,\ldots,v_j,\ldots,v_n)
\]
\end{prova}

\section{Algebraic and Geometric definitions}

What is a bit surprising is that these simple and natural properties imply the equivalence between the algebraic and geometric definitions of the determinants. Since $\det(e_1,\ldots,e_n) = 1$, it is enough to prove the following result.

\begin{proposicao}
$\det$ is multilinear and alternating.
\end{proposicao}

\begin{prova}
First observe that D1 and D2 imply the following property
\begin{enumerate}
\item[\textbf{D3}] For all vectors $v_1,\ldots,v_n$ and all $t$, we have that 
\[
\det(v_1,\ldots,v_i+tv_j,\ldots,v_j,\ldots,v_n) = \det(v_1,\ldots,v_i,\ldots,v_j,\ldots,v_n)
\]
\end{enumerate}
In fact, if $t$ is zero, this is trivially true and, if $t$ is not zero, we have that
\begin{eqnarray*}
\det(v_1,\ldots,v_i+tv_j,\ldots,v_j,\ldots,v_n)
& \stackrel{\text{D1}}{=} & t^{-1}\det(v_1,\ldots,v_i+tv_j,\ldots,tv_j,\ldots,v_n) \\
& \stackrel{\text{D2}}{=} & t^{-1}\det(v_1,\ldots,v_i,\ldots,tv_j,\ldots,v_n) \\
& \stackrel{\text{D1}}{=} & \det(v_1,\ldots,v_i,\ldots,v_j,\ldots,v_n)
\end{eqnarray*}
Second note that, if $\{v_1,\ldots,v_n\}$ is linearly dependent, then one of its elements can be written as linear combination of the others. Assuming that 
\[
v_i = t_1v_1 + \cdots + t_{i-1}v_{i-1} + t_{i+1}v_{i+1} + \cdots t_nv_n
\]
we have that
\begin{eqnarray*}
\det(v_1,\ldots,v_n) 
& = & \det(v_1,\ldots,v_i,\ldots,v_n) \\
& = & \det(v_1,\ldots,t_1v_1 + \cdots + t_{i-1}v_{i-1} + t_{i+1}v_{i+1} + \cdots t_nv_n,\ldots,v_n) \\
& \stackrel{\text{D3}}{=} & \det(v_1,\ldots,0,\ldots,v_n) \\
& \stackrel{\text{D1}}{=} & 0
\end{eqnarray*}
In order to prove the multilinearity, since we already have property D1, it is enough to show that
\[
\det(v_1,\ldots,v_i+v,\ldots,v_n) = \det(v_1,\ldots,v_i,\ldots,v_n) + \det(v_1,\ldots,v,\ldots,v_n)
\]
We need to analyse some possibilities. If all the sets
\[
\{v_1,\ldots,v_i+v,\ldots,v_n\}, \quad \{v_1,\ldots,v_i,\ldots,v_n\} \quad \mbox{and} \quad \{v_1,\ldots,v,\ldots,v_n\}
\]
are linearly dependent, it follows that
\[
\det(v_1,\ldots,v_i+v,\ldots,v_n) = 0 = \det(v_1,\ldots,v_i,\ldots,v_n) + \det(v_1,\ldots,v,\ldots,v_n)
\]
If $\{v_1,\ldots,v_i,\ldots,v_n\}$ é linearly independent, then it is a basis of $V$, and then we can write
\[
v = t_1v_1 + \cdots + t_iv_i + \cdots + t_nv_n
\]
Then it follows
\begin{eqnarray*}
\det(v_1,\ldots,v_i+v,\ldots,v_n)
& = & \det(v_1,\ldots,t_1v_1 + \cdots + (1+t_{i})v_{i} + \cdots t_nv_n,\ldots,v_n) \\
& \stackrel{\text{D3}}{=} & \det(v_1,\ldots,(1+t_{i})v_{i},\ldots,v_n) \\
& \stackrel{\text{D1}}{=} & (1+t_{i})\det(v_1,\ldots,v_{i},\ldots,v_n) \\
& \stackrel{\text{D1}}{=} & \det(v_1,\ldots,v_{i},\ldots,v_n) + \det(v_1,\ldots,t_iv_{i},\ldots,v_n) \\
& \stackrel{\text{D3}}{=} & \det(v_1,\ldots,v_{i},\ldots,v_n) + \det(v_1,\ldots,v,\ldots,v_n)
\end{eqnarray*}
For the remaining possibilities, we proceed as in the second one, with some minor adaptations. If $\{v_1,\ldots,v,\ldots,v_n\}$ is linearly independent, just change the roles of $v$ and $v_i$. If $\{v_1,\ldots,v_i+v,\ldots,v_n\}$ is linearly independent, write $v_i = v_i+v-v$ and change the roles of $v_i+v$ and $v_i$. In order to show the alternating property, it is enough to notice that
\begin{eqnarray*}
0
& = & \det(v_1,\ldots,v_i+v_j,\ldots,v_i+v_j,\ldots,v_n) \\
& = & \det(v_1,\ldots,v_i,\ldots,v_i,\ldots,v_n) + \det(v_1,\ldots,v_i,\ldots,v_j,\ldots,v_n) + \\
& & + \det(v_1,\ldots,v_j,\ldots,v_i,\ldots,v_n) + \det(v_1,\ldots,v_j,\ldots,v_j,\ldots,v_n) \\
& = & \det(v_1,\ldots,v_i,\ldots,v_j,\ldots,v_n) + \det(v_1,\ldots,v_j,\ldots,v_i,\ldots,v_n)
\end{eqnarray*}
where we used multilinearity in the second equality and also that
\[
\det(v_1,\ldots,v,\ldots,v,\ldots,v_n) = 0
\]
in the first and third inequality, since $\{v_1,\ldots,v,\ldots,v,\ldots,v_n\}$ is linearly dependent.
\end{prova}

From a didactic point of view, we can just introduce properties V1 and V2, justify them by presenting the respective pictures, introduce the determinants axiomatically through properties D1 and D2, prove the previous proposition, and then follow the standard approach.

\appendix

\section{Appendix}

\begin{proposicao}
Property V1 holds.
\end{proposicao}

\begin{prova}
The proof follows from the following decomposition
\[
P(v_1,\ldots,(s+t)v_i,\ldots,v_n) = P(v_1,\ldots,sv_i,\ldots,v_n) \mathop{\dot{\cup}} (sv_i + P(v_1,\ldots,tv_i,\ldots,v_n))
\]
since this implies that $f(t) = \vol(v_1,\ldots,tv_i,\ldots,v_n)$ is an additive and monotonous function and thus $f(t) = tf(1)$, for all nonnegative $t$, which is the property V1. The proof of the above decomposition is straightforward. 

In order to show the $\subset$ inclusion, let $v \in P(v_1,\ldots,(s+t)v_i,\ldots,v_n)$. Thus 
\[
v = t_1v_1 + \cdots + t_i(s+t)v_i + \cdots + t_nv_n
\]
where $t_1,\ldots,t_n \in [0,1)$. On one hand, if $t_i(s+t) < s$, we can write
\[
v = t_1v_1 + \cdots + \frac{t_i(s+t)}{s}sv_i + \cdots + t_nv_n
\]
showing that $v \in P(v_1,\ldots,sv_i,\ldots,v_n)$. On the other hand, if $t_i(s+t) \geq s$, we can write
\[
v = sv_ i + t_1v_1 + \cdots + \frac{t_i(s+t)-s}{t}tv_i + \cdots + t_nv_n
\]
showing that $v \in sv_i + P(v_1,\ldots,tv_i,\ldots,v_n)$, since in this case $\frac{t_i(s+t)-s}{t} \in [0,1)$. The other inclusion and the disjoint property in the above decomposition are easier to prove and left to the reader.
\end{prova}

\begin{proposicao}
Property V2 holds.
\end{proposicao}

\begin{prova}
Since $\mu$ is an additive measure invariant by translations, the proof follows immediately from the following decompositions
\[
P(v_1,\ldots,v_i,\ldots,v_j,\ldots,v_n) = A \mathop{\dot{\cup}} B
\]
and
\[
P(v_1,\ldots,v_i+v_j,\ldots,v_j,\ldots,v_n) = A \mathop{\dot{\cup}} (v_j + B)
\]
where
\[
A = \{ r_1v_1 + \cdots + r_i(v_i+v_j) + \cdots + r_jv_j + \cdots + r_nv_n : r_1,\ldots,r_n \in [0,1), r_i+r_j < 1\}
\]
and
\[
B = \{ s_1v_1 + \cdots + s_iv_i + \cdots + s_j(v_i+v_j) + \cdots + s_nv_n : s_1,\ldots,s_n \in [0,1), s_j < s_i+s_j < 1\}
\]

In order to show the $\subset$ inclusion in the first decomposition, we consider $v \in P(v_1,\ldots,v_i,\ldots,v_j,\ldots,v_n)$. Thus
\[
v = t_1v_1 + \cdots + t_iv_i + \cdots + t_jv_j + \cdots + t_nv_n
\]
where $t_1,\ldots,t_n \in [0,1)$. On one hand, if $t_j \geq t_i$, we can write 
\[
v = t_1v_1 + \cdots + t_i(v_i+v_j) + \cdots + (t_j-t_i)v_j + \cdots + t_nv_n \in A
\]
since $t_i + t_j - t_i < 1$. On the other hand, if $t_j < t_i$, we can write
\[
v = t_1v_1 + \cdots + (t_i-t_j)v_i + \cdots + t_j(v_i+v_j) + \cdots + t_nv_n \in B
\]
since $t_j < t_i - t_j + t_j < 1$. The other inclusion and the disjoint property in the first decomposition are easier to prove and left to the reader. 

In order to show the $\subset$ inclusion in the second decomposition, we consider $v \in P(v_1,\ldots,v_i+v_j,\ldots,v_j,\ldots,v_n)$. Thus
\[
v = t_1v_1 + \cdots + t_i(v_i+v_j) + \cdots + t_jv_j + \cdots + t_nv_n
\]
where $t_1,\ldots,t_n \in [0,1)$. On one hand, if $t_i+t_j < 1$, we have that $v \in A$. On the other hand, if $t_i+t_j \geq 1$, we can write
\[
v = v_j + t_1v_1 + \cdots + (1-t_j)v_i + \cdots + (t_i+t_j-1)(v_i+v_j) + \cdots + t_nv_n \in B
\]
since $t_i+t_j-1 < 1-t_j+t_i+t_j-1 < 1$. The other inclusion and the disjoint property in the second decomposition are easier to prove and left to the reader. 
\end{prova}

\end{document}